\documentclass{article}

\usepackage{enumerate}
\usepackage{amsmath,amssymb}
\usepackage{ascmac}
\usepackage{amscd}

\usepackage[top=30truemm,bottom=30truemm,left=40truemm,right=40truemm]{geometry}

\newtheorem{theorem}{Theorem}[section]
\newtheorem{Cor}{Corollary}[section]
\newtheorem{Prop}{Proposition}[section]
\newtheorem{Def}{Definition}[section]

\newtheorem{lemma}{Lemma}[section]

\newtheorem{Rem}{Remark}[section]
\newtheorem{proof}{Proof.}

\numberwithin{equation}{section}

\title{Generalized $k$-regular sequences III:\\ Arithmetical properties of generalized $k$-regular series }

\author{Eiji Miyanohara \\
E-mail: j1o9t5acrmo@fuji.waseda.jp}
\date{\today}

\begin{document}

\date{}
\maketitle

\begin{abstract}
Let $F(z)$ be a $k$-regular series in $\mathbb{Z}[[z]]$ and $b$ be an integer with $b\ge2$.
Bell, Bugeaud and Coons [BelBC] proved that $F(\frac{1}{b})$ is either rational or transcendental.
In [Mi], we introduce a generalized $k$-regular sequence as a unification of several kinds of important sequences including $k$-regular, $k$-additive and $k$-multiplicative sequences.
In this paper, we give a  generalization of the result of Bell, Bugeaud and Coons for certain generalized $k$-regular series.
Especially, we show that the values of irrational generating functions of certain sum of $k$-additive sequences and certain $k$-multiplicative sequences are either rational or transcendental.
Moreover, we also give a partly generalization of a result obtained by Tachiya[Ta].
Especially, we show that the values of irrational generating functions of certain $k$-additive sequences and certain $k$-multiplicative sequences give transcendental numbers.
\end{abstract}

\section{Introduction}
\quad
Let ${\bf{a}}:={(a(n ))}_{n\ge0}$ be a given sequence.
For any non-negative integer $e$, set $$S_e({\bf{a}}):=\{{(a(k^en+j ))}_{n\ge0}\; |0\le j \le k^e-1 \}.$$
 Allouche and Shallit [AlS] introduced the notion of $k$-regular sequence as follows. A sequence ${(a(n ))}_{n\ge0}$ is defined to be $k$-regular if the set $S$ is contained in a finitely genelated $\mathbb{Q}$-module of sequences.
Allouche-Shallit also proved that the set of generating function of $k$-regular sequence (called a $k$-regular series) forms a ring under the usual addition and the canonical convolution.
Later Becker [Bec] and Nishioka [Ni], which characterizes $k$-regular sequences by using the $k$-regular series. 
(See Theorem $5.1.2$ in [Ni].)
\begin{theorem}\label{inf}{\rm[Bec,Ni]}
A sequence ${(a(n ))}_{n\ge0}$ is $k$-regular if and
only if there exist a positive integer $d$, $d$ power sereis $f_1(z)\cdots f_d(z)\in \mathbb{Q}[[z]]$ with
$f_1(z) = f(z)$ given in and a $d\times d$ matrix $A(z)$ whose entries are polynomials
in $z$ of degrees less than $k$ with coefficients in $\mathbb{Q}$ such that

\begin{align}\label{2}
\left(
    \begin{array}{c}
      f_1(z) \\
      f_2(z) \\
      \vdots \\
      f_{d}(z)
    \end{array}
  \right)=A(z)
  \left(
    \begin{array}{c}
      f_1(z^k) \\
      f_2(z^k)\\
      \vdots \\
      f_{d}(z^k)
   \end{array}
  \right).
  \end{align}
\end{theorem}
By Theorem \ref{inf}, the $k$-regular series can be regard as Mahler function. (See chapter $5$ in [Ni]).
Therefore, the arithmetical properties  of the special value of $k$-regular series was investigated in Mahler function theory.
Recently, Bell, Bugeaud and Coons [BelBC] proved the following theorem. (See Theorem $8.1$ in [BelBC] or Theorem $2.5.1$ in [CoS].)
\begin{theorem}\label{BelBC}{\rm[BelBC]}
Let $F(z)$ be a $k$-regular series in $\mathbb{Z}[[z]]$ and $b$ be an integer with $b\ge2$.
Then $F(\frac{1}{b})$ is either rational or transcendental.
\end{theorem}
The proof of Theorem \ref{BelBC} relies on $p$-adic Schmidt subspace theorem.\\
\quad
On the other hands, Gel'fond [Gel] introduced the two functions related with the base $k$-representation as follows.
A sequence ${(a(n))}_{n\ge0}$ is $k$-additive if and only if, for any non-negative integers $e$, $n$ and $j$ with $0\le j\le k^e-1$, ${(a(n))}_{n\ge0}$ satisfies the following additive relation
\begin{align}\label{4}a(k^en+j )=a(k^en)+a(j)
\end{align}
and $a(0)=0$.
A sequence ${(a(n))}_{n\ge0}$ is $k$-multiplicative if and only if, for any non-negative integers $e$, $n$ and $j$ with $0\le j\le k^e-1$, ${(a(n))}_{n\ge0}$ satisfies the following multiplicative relation
\begin{align}\label{5}a(k^en+j )=a(k^en)a(j)
\end{align}
and $a(0)=1$.\\
Recently, we introduce a generalized $k$-regular sequence as a unification of several kinds of important sequences including $k$-regular, $k$-additive and $k$-multiplicative sequences in [Mi]. 
\begin{Def}{\rm[Mi]}\label{defofk-pro}
 A sequnece ${(a(n ))}_{n\ge0}$ is {\it generalized $k$-regular} if and only if, there exist an integer $d$, for any non-negative integer $e$, the dimention of the $\mathbb{Q}$-module of a sequence generated by $S_e$ is at most $d$. The generating series of generalized $k$-regular sequences is called
a {\it generalized $k$-regular series}. To show the role of $d$ more preciesly, a generalized $k$-regular sequence or series is also called
a {\it generalized $(k,d)$-regular sequence or series}, respectively.
\end{Def}
In [Mi], we give the following generalization of Theorem \ref{inf} for generalized $k$-regular sequences as follows.
\begin{theorem}{\rm[Mi]}\label{th2}
 A sequence ${(a(n ))}_{n\ge0}$ is generalized ($k$, $d$)-regular if and
only if, for any non-negative integer $e$, there exist a positive integer $d$, $d$ power sereis $f_{e,1}(z),f_{e,2}(z),$
$\cdots,f_{e,d}(z)\in \mathbb{Q}[[z]]$ with
$f_{0,1}(z)$ being $f(z)$ given in and a $d\times d$ matrix $A_e(z)$ whose entries are polynomials
in $z$ of degrees less than $k$ with coefficients in $\mathbb{Q}$ such that
\begin{align}\label{8}
&\left(
    \begin{array}{c}
      f_{e,1}(z) \\
      f_{e,2}(z) \\
      \vdots \\
      f_{e,d}(z)
    \end{array}
  \right)=A_e(z)\left(
    \begin{array}{c}
      f_{e+1,1}(z^k) \\
      f_{e+1,2}(z^k)\\
      \vdots \\
      f_{e+1,d}(z^k)
   \end{array}
  \right) \;\;(e\ge0).
  \end{align} 
\end{theorem}
Now we give the natural three examples of Theorem \ref{th2}.
Let $g_0(z)$ be the generating function of a $k$-multiplicative sequence ${(a(n ))}_{n\ge0}$ (generalized ($k$,$1$)-regular series).
For any non-negative integer $e$, we define $g_{e}(z)$ as $g_{e}(z):=\sum_{n=0}^\infty a(k^en)z^n$.
The series $g_0(z)$ has the following infinite chains equations
\begin{align}\label{6}
g_e(z)=(\sum_{j=0}^{k-1} a(jk^e)z^{j})g_{e+1}(z^k). 
\end{align}
The arithmetical properties of infinite product \eqref{6} was investigated in Mahler function theory. (See [AmV1,AmV2,Ta].)
Let $h_0(z)$ be the generating function of $k$-additive sequence ${(b(n))}_{n\ge0}$ (generalized ($k$,$2$)-regular series).
For any non-negative integer $e$, we define $h_{e}(z)$ as $h_{e}(z):=\sum_{n=0}^\infty b(k^en)z^n$. The series $h_{0}(z)$ has the following infinite chains matrix equations
\begin{align}\label{111111111}
    &\left(\begin{array}{c}
      h_e(z) \\
      \frac{1}{1-z}
      \end{array}
  \right) =\begin{pmatrix}
\sum_{j=0}^{k-1} z^j&\sum_{j=0}^{k-1} b(jk^e)z^j \\
0 &\sum_{j=0}^{k-1} z^j
\end{pmatrix}\left(\begin{array}{c}
      h_{e+1}(z^k) \\
      \frac{1}{1-z^k}
      \end{array}
  \right).
\end{align}
The set of generalized $k$-regular sequences forms a ring under the usual
addition and the canonical convolution. (See Theorem 2.2 in [Mi].) Therefore, the power series $g_0(z)+h_0(z)$ is also a generalized $k$-regular series. The power series $g_0(z)+h_{0}(z)$ has the following infinite chains matrix equations
{\footnotesize{\begin{align}\label{1111111111}
    &\left(\begin{array}{c}
      g_e(z)+h_e(z) \\
      g_e(z)
      \\
      h_e(z)\\
      \frac{1}{1-z}
      \end{array}
  \right) =\begin{pmatrix}
\sum_{j=0}^{k-1} (a(jk^e)+1)z^{j}&-\sum_{j=0}^{k-1} z^j&-\sum_{j=0}^{k-1} a(jk^e)z^{j}
& \sum_{j=0}^{k-1} b(jk^e)z^j \\
0&\sum_{j=0}^{k-1} a(jk^e)z^{j}&0&0\\
0&0 &\sum_{j=0}^{k-1} z^j&\sum_{j=0}^{k-1} b(jk^e)z^j\\
0&0 &0&\sum_{j=0}^{k-1} z^j
\end{pmatrix}\left(\begin{array}{c}
      g_{e+1}(z^k)+h_{e+1}(z^k) \\
      g_{e+1}(z^k)
      \\
      h_{e+1}(z^k)\\
      \frac{1}{1-z^k}
      \end{array}
  \right).
\end{align}}}
\\
\quad
The purpose of this paper investigates the arithmetical properties of certain generalized $k$-regular series as follows.\\
\quad
We denote a $(i,j)$-componet of $A_e(z)$ by $\sum_{s=0}^{k-1} a_{e,s,i,j} z^s$.
We assume that ,for any non-negative integers $e$ and $j$ with $1\le j\le d$, there exists a positive constant $C$
\begin{align}\label{asy}
\;|f_{e,j}(0)|\le C
\end{align}
and, for any $\epsilon>0$, there exists an integer $N(\epsilon)$ such that, for any $e\ge N(\epsilon)$, $i,j$ with $1\le i, j\le d$ and $s$ with $0\le s\le k-1$,
\begin{align}\label{asy21}
|a_{e,s,i,j}|\le {e}^{\epsilon k^e}.
\end{align}
Moreover, we assume that, for any non-negative integer $e$,
\begin{align}\label{ininteger}
A_e(z)\in{\mathbb{Z}[z]}^{d\times d}.
\end{align}
\begin{theorem}\label{main2}
Let $b$ be an integer with $b\ge 2$ and $f(z)=f_{0,1}(z)$ be satisfies the equations \eqref{8} with \eqref{asy}, \eqref{asy21} and \eqref{ininteger}. Then $f(\frac{1}{b})$ is either rational or transcendental.
\end{theorem}
We prove Theorem \ref{main2} by modifying the method of proof of Theorem \ref{BelBC}. (See the proof of Theorem $2.5.1$ in [CoS].)
By \eqref{1111111111}, we get the following corollary of Theorem \ref{main2}.
\begin{Cor}
Let $g_0(z)$ and $h_0(z)$ be defined by the above and $b$ be an integer with $b\ge 2$.
Assmue that, for any non-negative integers $e$ and $j$ with $0\le j\le k-1$, $a(k^e)$ and $b(k^e)$ are an integers, $a(jk^e)$ and $b(jk^e)$ satisfy \eqref{asy21}. Then $g_0(\frac{1}{b})+h_0(\frac{1}{b})$ is either rational or transcendental.
\end{Cor}
By Theorem \ref{main2} and the most classical Mahler method (See 20p in [Ma]),
we also prove the following theorem.
\begin{theorem}\label{main}
Let $b$ be an integer with $b\ge 2$ and irrational powers series $f(z)=f_{0,1}(z)$ be satisfies the equations \eqref{8} with \eqref{asy}, \eqref{asy21} and \eqref{ininteger}. Assume that, for any non-negative integer $e$, $\det A_e(\frac{1}{b^{k^e}})\neq 0$.
Then at least one among the numbers $f(\frac{1}{b})=f_{0,1}(\frac{1}{b}),f_{0,2}(\frac{1}{b}),$ $\cdots, f_{0,d}(\frac{1}{b})$ is transcendental.
\end{theorem}
Theorem \ref{main} gives a partly generalization of a result obtained by Theorem $1$ in [Ta].
By \eqref{6}, we get the following corollary of Theorem \ref{main}.
\begin{Cor}Let $g_0(z)$ be defined by the above with irrational and $b$ be an integer with $b\ge 2$. Assmue that, for any non-negative integers $e$ and $j$ with $0\le j\le k-1$, $a(k^e)$ is an integer, $a(jk^e)$ satisfies \eqref{asy21} and $\sum_{j=0}^{k-1} a(jk^e)\frac{1}{b^{jk^e}}\neq 0$. Then $g_0(\frac{1}{b^{}})$ is transcendental.
\end{Cor}
This corollary is covered by Theorem $1$ in [Ta].
By \eqref{111111111}, we get the following corollary of Theorem \ref{main}.
\begin{Cor}
Let $h_0(z)$ be defined by the above with irrational and $b$ be an integer with $b\ge 2$. Assmue that, for any non-negative integers $e$ and $j$ with $0\le j\le k-1$, $b(k^e)$ is an integer and $b(jk^e)$ satifies \eqref{asy21}.
Then $h_0(\frac{1}{b^{}})$ is transcendental.
\end{Cor}
This corollary is new.\\
\quad
This paper is organized as follows.
In section $2$, we gather lemmas for the proof of the theorems and the proposition.
In section $3$, we give a proof of Theorem \ref{main2}.
In section $4$, we give a proof of Theorem \ref{main}.
In section $5$, we give the other examples of Theorem \ref{main} by related with the certain digital pattern sequences.
\section{Preliminaries}
In this section, we gather lemmas for the proof of the theorems.
The following lemma is need for the proof of Theorem \ref{main2}.
The following lemma is known as Siegel's lemma. (See Lemma $1.4.2$ in [Ni].)
\begin{lemma}[Siegel's lemma]\label{Siegel}
Consider the $m$ equations in $n$ unknowns
\begin{align}\label{21}
a_{k1}x_1+\cdots+a_{kn}x_n=0\;\;k = 1,2,\cdots,m
\end{align}
with rational integral coefficients $a_{ij}$, and with $0 < m < n$. Let $A$ be
a positive integer such that $A\ge |a_{ij}|$, for all $i$ and $j$. Then there is a
nontrivial solution $x_1, x_2,\cdots , x_n $in rational integers of equations \eqref{21}
such that
\begin{align}\label{22}
|x_j|<1+(nA)^{n/(n-m)}\;\;j = 1,2,\cdots,n.
\end{align}
\end{lemma}
The following lemma is need for the proof of Theorem \ref{main2}.
The following lemma is known as $p$-adic Schmidt subspace theorem. (See Theorem $E.10$ in [Bu] or Theorem $2.5.4$ in [CoS].)
\begin{lemma}[$p$-adic Schmidt subspace theorem ]\label{Subspcace}
Let $n\ge 2$, $\delta>0$. and
let $p_1,\cdots,p_s$ be distinct prime numbers. Further, let $L_{1\infty},\cdots,L_{n\infty}$ be linearly
independent linear forms in $X_1,\cdots,X_n$ with algebraic coefficients in $\mathbb{C}$, and for
$j=1,\cdots, s$, $L_{1j},\cdots,L_{nj}$ be linearly
independent linear forms in $X_1,\cdots,X_n$ with algebraic coefficients in ${\bar{\mathbb{Q}}}_p$.
Consider the inequality
\begin{align}\label{23}
|L_{1\infty}({\bf{x}})\cdots L_{n\infty}({\bf{x}})|\prod_{j=1}^s |L_{1j}({\bf{x}})\cdots L_{nj}({\bf{x}})|_p<|\max\{x_1,\cdots, x_m\}|^{-\delta}
\end{align}
with ${\bf{x}}:=(x_1,\cdots, x_m)$ in ${\mathbb{Z}}^n$.
There are a finite number of proper linear subspaces $T_1,\cdots T_t$ of ${\mathbb{Q}}^n$
such that all solutions of \eqref{23} lie in $T_1\cup\cdots\cup T_t$.
\end{lemma}
The following notion is need for the construct of the examples of Theorem \ref{main}.
(See Definition $1$ in [AmV1].)
\begin{Def}{\rm[AmV1]}\label{irrmea}
Let $f(z)\in K[[z]]$. We define the irrationality measure $\mu(f)$ to be the infimum of $\mu$ such that;
\begin{align*}
ord(A(z)f(z)-B(z))\le \mu M
\end{align*}
for all nonzero $A(z), B(z)\in K[z]$ with $\max(degA(z), degB(z))\le M$ (for $M\ge M_0$, some $M_0$ depend only on $f(z)$). If there does not exist such a $\mu$, $\mu(f):=+\infty$.
\end{Def}
The following lemma is need for the construct of the examples of Theorem \ref{main}.(See Theorem $5$ in [DuN].) 
\begin{lemma}{\rm[DuN]}\label{DN}
Let $K$ be a commutative field and $c_1$, $c_2$, $c_3$ be real numbers with
$0 <c_1<c_2$, $c_3 \ge1$. Let ${(m(n))}_{n\ge0}$ be an increasing sequence of nonnegative integers satisfying
$m(n+1)-m(n)\le c_3$. Let $k\ge2$ be an integer and $f(z)\in K[[z]]$ . Suppose that for large positive integer $n$ there
exists a sequence ${(P_n(z),Q_n(z))}_{n=0}^\infty$ in $ {K[z]}^2$ satisfying
\begin{align*}
&P_n(z)Q_{n+1}(z) -P_{n+1}(z)Q_n(z) \neq 0,\\
&degQ_n(z), degP_n(z)\le c_1k^{m(n)},\\
&ord(Q_n(z)f (z)-P_n(z))\ge c_2k^{m(n)}. 
\end{align*}
Then $\mu(f)<+\infty$.
\end{lemma}
\section{Proof of Theorem \ref{main2}}\label{s:3}
In this section, we prove Theorem \ref{main2}.
Let $p$ be a positive integer parameter with $p>d+5$.
We shall denote by $c_1, c_2, \cdots$ positive constants independent of $\epsilon,p,e$.
For any non-negative integers $e$ and $j$ with $1\le j\le d$, we define the $(a_{e,j}(n))_{n\ge 0}$ by $f_{e,j}(z)=\sum_{\infty}a_{e,j}(n)z^n$.
\begin{lemma}\label{31}
Notation is the same as for section $1$.
Then, for any $e\ge N(\epsilon)$, $j$ with $1\le j\le d$ and $n\ge0$,
\begin{align}\label{asy2}
|a_{e,j}(n)|\le  {e}^{\epsilon k^e(1+n)}.
\end{align}
\end{lemma}
\begin{proof}{\rm
By \eqref{asy} and \eqref{asy21}, one can show analogously to the proof of Lemma 3 in [Ta].}
\end{proof}
\begin{lemma}\label{siegel2}
Notation is the same as for section $1$.
For any $e\ge N(\epsilon)$, $j$ with $1\le j\le d$ and $n\ge0$
there exist auxiliary functions 
for any $e\ge N(\epsilon)$, $j$ with $1\le j\le d$ and $n\ge0$, we have
\begin{align}\label{aux}
Q_e(z)f_{e,j}(z)-P_{e,j}(z)=z^{dp+p+1}G_{e,j}(z).
\end{align}
with polynomials $Q_e (z)=\sum_{i=0}^{dp} q_{e,i} z^i, P_{e,j}(z)=\sum_{i=0}^{dp} p_{e,j,i} z^i,\in\mathbb{Z}[z]/0$ of degrees at most $dp$ and $G_{e,j}(z)=\sum_{n=0}^{\infty}g_{e,j}(n)z^n$, such that
\begin{align}\label{aux2}
&|q_{e,i}|\le 1+((dp+1){e}^{\epsilon k^e(1+dp)})^{(dp+1)}\le{e}^{\epsilon c_1p^2 k^e},\\
&|p_{e,j,i}|\le (dp+1) (1+((dp+1){e}^{\epsilon k^e(1+dp)})^{dp}){e}^{\epsilon k^e(1+dp)}\le{e}^{\epsilon c_2p^2 k^e},\\
& |g_{e,j}(n)|\le (dp+1) ({1+((dp+1){e}^{\epsilon k^e(1+dp)})^{dp}}{e}^{\epsilon k^e(1+n)})\le{e}^{\epsilon (c_3p^2+n) k^e}.
\end{align}
\end{lemma}
\begin{proof}{\rm
By Lemma \ref{21} and \ref{31}, one can show analogously to the proof of Lemma 5 in [Am2].}
\end{proof}
\begin{lemma}\label{siegel3}
Notation is the same as for section Lemma \ref{siegel2}.
Then, for any $e\ge N(\epsilon)$, $j$ with $1\le j\le d$,
\begin{align}\label{aux10}
ord \;Q_e(z)\le ord \; P_{e,j}(z).
\end{align}
\end{lemma}
\begin{proof}{\rm
We denote $ord \;Q_e(z)$ by $J_e$ and $Q^{'}_e(z)\in \mathbb{Z}[z]$ by $Q_e(z)=z^{J_e}Q^{'}_e(z)$.
\begin{align}\label{aux3}
Q^{'}_e(z)f_{e,j}(z)-\frac{P_{e,j}(z)}{z^{J_e}}=z^{dp+p+1-J_e}G_{e,j}(z).
\end{align}
By the definition of $J_e$, 
\begin{align}\label{aux4}
dp+p+1-J_e\ge0.
\end{align}
By \eqref{aux4} and right hands side of \eqref{aux3}, $\frac{P_{e,j}(z)}{z^{J_e}}\in \mathbb{Z}[z]$.
Thereofere, we get \eqref{aux10}.}
\end{proof}
There exists an integer $J$ with $J\le dp$ such that
\begin{align}\label{J}
\#\{e\;|\; ord \;Q_e(z)=J\}=\infty.
\end{align}
We denote the set $\{e\;|\; ord \;Q_e(z)=J\}$ by $B$. By \eqref{aux}, for any integer $e$ in $B$ and $j$ with $1\le j\le d$, we have
\begin{align}\label{aux11}
Q^{'}_e(z)f_{e,j}(z)-\frac{P_{e,j}(z)}{z^{J}}=z^{dp+p+1-J}G_{e,j}(z).
\end{align}
We replace $Q^{'}_e(z)$ and $\frac{P_{e,j}(z)}{z^{J}}\in\mathbb{Z}[z]$ by $Q_e(z)$ and $P_{e,j}(z)$.
By \eqref{8} and \eqref{aux11}, there exist polynomials $a_{e,j,0}(z)$ ($1\le j\le d$) with degrees at most $k^e$ such that
\begin{align}\label{aux6}
Q_e(z^{k^e})f(z)-\sum_{j=1}^d a_{e,j,0}(z)P_{e,j}(z^{k^e})=z^{(dp+p+1-J){k^e}}\sum_{j=1}^d a_{e,j,0}(z)G_{e,j}(z^{k^e}).
\end{align}
\begin{lemma}\label{siegel5}
If $|z|\le\frac{2}{3{e}^{\epsilon c_1p^2}}$, then, for sufficiently large integer $e$, 
\begin{align}\label{pl}
|Q_e(z^{k^e})|\ge 1/2.
\end{align}.
\end{lemma}
\begin{proof}{\rm
By the definition of $Q_e(z)$ and \eqref{aux2}, we have
\begin{align}\label{lower}
|Q_e(z^{k^e})|\ge 1-\frac{{(2e)}^{\epsilon c_1p^2 k^e}}{{(3e)}^{\epsilon c_1p^2k^e}}-\cdots-\frac{{(2e)}^{\epsilon c_1p^2 k^e}}{{(3e)}^{\epsilon c_1p^2(dp-J)}}\ge\frac{1}{2}.
\end{align}}
\end{proof}
\begin{lemma}\label{siegel4}
Let $b$ be an integer with $b\ge 2$. If 
\begin{align}\label{aux81}\frac{3{e}^{2\epsilon c_1dp^2}}{2}<{b}\end{align}
then, for sufficiently large integer $e$,
\begin{align}\label{aux8}
|f(\frac{1}{b})-\frac{\sum_{j=1}^d a_{e,j,0}(\frac{1}{b^{}})P_{e,j}(\frac{1}{b^{k^e}})}{Q_e(\frac{1}{b^{k^e}})}|\le \frac{1}{b^{(dp+p+1-J){k^e}}}
2C_1(\epsilon)d^e {e}^{\epsilon k^e}\frac{b}{b-1}  {e}^{\epsilon c_3p^2}\frac{{b}^{k^e}}{{b}^{k^e}-{e}^{\epsilon k^e}},
\end{align}
where $C_1(\epsilon)$ is a positive constant independent of $e$.
In paticular, 
\begin{align}\label{aux9}\lim_{e\to \infty} \frac{\sum_{j=1}^d a_{e,j,0}(\frac{1}{b})P_{e,j}(\frac{1}{b^{k^e}})}{Q_e(\frac{1}{b^{k^e}})}=\lim_{e\to \infty} \frac{\sum_{j=1}^d a_{e,j,0}(\frac{1}{b})P_{e,j}(\frac{1}{b^{k^e}})}{q_{e,0}}=f(\frac{1}{b}).
\end{align}
\end{lemma}
\begin{proof}{\rm
By $q_{e,0}\neq 0$, \eqref{aux6} and \eqref{aux2}, for any sufficiently large $e$, we have
\begin{align}\label{con}
&|Q_e(\frac{1}{b^{k^e}})f(\frac{1}{b})-\sum_{j=1}^d a_{e,j,0}(\frac{1}{b})P_{e,j}(\frac{1}{b^{k^e}})|\le\frac{1}{b^{(dp+p+1-J){k^e}}}
|\sum_{j=1}^d a_{e,j,0}(\frac{1}{b})|\sum_{n=0}^\infty  \frac{{e}^{\epsilon (c34p^2+n) k^e}}{{b}^{nk^e}}\nonumber\\
&\le \frac{1}{b^{(dp+p+1-J){k^e}}}
C_1(\epsilon)d^e{e}^{\epsilon k^e}\frac{b}{b-1} \sum_{n=0}^\infty  \frac{{e}^{\epsilon (c_3p^2+n) k^e}}{{b}^{nk^e}}\le
\frac{1}{b^{(dp+p+1-J){k^e}}}
C_1(\epsilon)d^e{e}^{\epsilon k^e}\frac{b}{b-1}  {e}^{\epsilon c_3p^2}\frac{{b}^{k^e}}{{b}^{k^e}-{e}^{\epsilon k^e}}.
\end{align}
By Lemma \ref{pl},\eqref{aux81} and \eqref{con}, we get
\begin{align}\label{con2}
&|f(\frac{1}{b^{k^e}})-\frac{\sum_{j=1}^d a_{e,j,0}(\frac{1}{b^{}})P_{e,j}(\frac{1}{b^{k^e}})}{Q_e(\frac{1}{b^{k^e}})}|\le
\frac{1}{b^{(dp+p+1-J){k^e}}}
2C_1(\epsilon)d^e{e}^{\epsilon k^e}\frac{b}{b-1}  {e}^{\epsilon c_3p^2}\frac{{b}^{k^e}}{{b}^{k^e}-{e}^{\epsilon k^e}}.
\end{align}
Moreover,by $(3.3)$, we have
\begin{align}\label{con3}
|Q_e(\frac{1}{b^{k^e}})-q_{e,0}|\le dp (\frac{{e}^{\epsilon c_4p^2}}{b})^{k^e}.
\end{align}
By \eqref{aux81},\eqref{con2} and \eqref{con3}, we get \eqref{aux9}.}
\end{proof}
Assume that \eqref{aux81}.
For any non-negative integer $e$, we define the integer tuples $(D_{e,0},\cdots,$ $D_{e,dp-J},D_{e,dp-J+1})$ as follows
\begin{align}
(D_{e,0},\cdots,D_{e,dp-J},D_{e,dp-J+1}):=(b^{(dp+1-J)k^{e}}q_{e,0},\cdots,b^{k^{e}}q_{e,dp-J},b^{(dp+1-J)k^{e}}\sum_{j=1}^d a_{e,j,0}(\frac{1}{b^{}})P_{e,j}(\frac{1}{b^{k^e}})).
\end{align}
From $(3.3)$ and $(3.4)$, we have
\begin{align}
\max\{D_{e,0},\cdots,D_{e,dp-J},D_{e,dp-J+1}\}\le b^{(dp+1-J)k^{e}} {e}^{\epsilon c_2p^2 k^e} 2C_1(\epsilon)d^e {e}^{\epsilon k^e}\frac{b}{b-1}\le b^{(dp+3)k^{e}}.
\end{align}
Moreover, there exist the integer sets $T:=\{s_1,s_2,\cdots,s_l\}$ with $0\le s_1<s_2<\cdots<s_l\le dp-J+1$ such that
\begin{align}
\#\{e\;|\; D_{e,i}\neq 0\;\mbox{for  $i$ in $T$ and}\; D_{e,i}=0\;\mbox{for  $i$ in $\{0,\cdots,dp-J+1\}/T$}\}=\infty
\end{align}
We put the set $E:=\{e\;|\; D_{e,i}\neq 0\;\mbox{for  $i$ in $T$ and}\; D_{e,i}=0\;\mbox{for  $i$ in $\{0,\cdots,dp-J+1\}/T$}\}$. We assume that $f(\frac{1}{b})\neq 0$ is an algebraic number. By $f(\frac{1}{b})\neq 0$, $q_{e,0}\neq 0$ and \eqref{aux9}, we have
\begin{align}\label{vv}
s_1=0,s_l=dp-J+1.
\end{align}
Let $S$ be the set of prime factor of $b$.
We define the linear form
\begin{align}\label{linear1}
L_{i,\infty}=x_i\;\; (1\le i\le l-1)
\end{align}
and 
\begin{align}\label{linear2}
L_{l,\infty}=f(\frac{1}{b^{}})\sum_{i=1}^{l-1} x_i+x_l.
\end{align}
Moreover, for any prime $p$ in $S$, we define the linear form
\begin{align}\label{linear3}
L_{i,p}=x_i\;\; (1\le i\le l).
\end{align}
For any sufficiently large integer $e$ in $E$, we define $(x_1,\cdots,x_l):=(D_{e,0},D_{e,s_2},\cdots,D_{e,dp-J+1})$.
\begin{align}\label{linear4}
&|L_{1\infty}({\bf{x}})\cdots L_{l\infty}({\bf{x}})|\prod_{p\in S}\prod_{j=1}^l |L_{1j}({\bf{x}})\cdots L_{lj}({\bf{x}})|_p\le
\frac{1}{b^{p{k^e}}}
2C_1(\epsilon)d^e{e}^{\epsilon k^e}\frac{b}{b-1}  {e}^{\epsilon c_4p^2} {e}^{\epsilon c_3dp^3 k^e}\frac{{b}^{k^e}}{{b}^{k^e}-{e}^{\epsilon k^e}}\nonumber\\&\le \frac{1}{b^{(p-3)k^e}}\le
{(\frac{1}{b^{(dp+3)k^{e}}})}^{\frac{p-3}{dp+3}}\le \frac{1}{{(\max\{D_{e,0},\cdots,D_{e,dp-J},D_{e,dp-J+1}\})}^{\frac{p-3}{dp+3}}}.
\end{align}
By \eqref{linear4} and Lemma \ref{Subspcace}, $(x_1,\cdots,x_l):=(D_{e,0},D_{e,s_2},\cdots,D_{e,dp-J+1})$ with in     $E$ lie in finitely many proper linear subspaces of ${\mathbb{Q}}^l$.
There exist an infinite set of distinct positive integers $E^{'}\subset E$ and a nonzero integer triple $(z_1, \cdots, z_l)$
such that
\begin{align}\label{linear5}
&z_1b^{(dp+1-J)k^{e}}q_{e,0}+z_2D_{e,s_2}+\cdots+z_lb^{(dp+1-J)k^{e}}\sum_{j=1}^d a_{e,j,0}(\frac{1}{b})P_{e,j}(\frac{1}{b^{k^e}})\nonumber\\&=z_1D_{e,0}+z_2D_{e,s_2}+\cdots+z_lD_{e,dp-J+1}=0, \;\; \mbox{for any}\; e\in E^{'}.
\end{align}
We define the integer $m$ as $m:=\min\{ i \;|z_i\neq 0\}$.
If $1<m$, we have
\begin{align}\label{linear6}
&z_mb^{(dp+1-J-s_m)k^{e}}q_{e,s_m}+\cdots+z_{l-1}b^{s_{l-1}k^{e}}q_{e,s_{l-1}}=z_lb^{(dp+1-J)k^{e}}\sum_{j=1}^d a_{e,j,0}(\frac{1}{b^{}})P_{e,j}(\frac{1}{b^{k^e}}), \;\; \mbox{for any}\; e\in E^{'}.
\end{align}
By\eqref{linear6}, $|q_{e,0}|\ge1$ and \eqref{aux9}, for sufficiently large $e$, we have
\begin{align}\label{limm}
|\sum_{j=1}^d a_{e,j,0}(\frac{1}{b^{}})P_{e,j}(\frac{1}{b^{k^e}})|\ge \frac{|f(\frac{1}{b})|}{2}\neq 0.
\end{align}
By \eqref{linear6} and \eqref{limm}, $e$ tend to infinity, we get
\begin{align}\label{linear10}
z_l=0.
\end{align}
By \eqref{linear6} and \eqref{linear10}
\begin{align}\label{linear11}
z_m=0.
\end{align}
This contradicts the definition of $m$.
Therefore, $m=1$.
Dividing \eqref{linear5} by $b^{(dp+1-J)k^{e}}q_{e,0}$ and \eqref{aux9}, $e$ tend to infinity, we get
\begin{align}\label{linear12}
&z_1+z_lf(\frac{1}{b})=0.
\end{align}
By \eqref{linear12}, $f(\frac{1}{b})$ is a rational number. This completes the proof of Theorem \ref{main2}.\hfill
$\square$
\section{Proof of Theorem \ref{main} }
Now we prove Theorem \ref{main}.
For any non-negative integer $e$, we define the matrix $B_e(z)$ as $B_e(z)=A_0(z)\cdots A_{e-1}(z^{k^{e-1}})$.
\begin{lemma}\label{invert}
For any sufficiently large integer $e$,
there exists an non-zero integer $D_e$ such that $D_eB^{-1}_e(\frac{1}{b})\in\mathbb{Z}^{d\times d}$ and 
\begin{align}\label{upper}D_e\le C_2(\epsilon)d!^e b^{dk^e} e^{\epsilon dk^e}
\end{align}
where $C_2(\epsilon)$ is a positive constant independent of $e$.
\end{lemma}
\begin{proof}{\rm By the computation of numerator of $\det A_i(z)$ with \eqref{asy} and \eqref{asy21}.}
\end{proof}
We assume that
\begin{align}\label{asss}
b>e^{\epsilon c_3p^2}
\end{align}
and $f(\frac{1}{b})=f_{0,1}(\frac{1}{b}),f_{0,2}(\frac{1}{b}),\cdots, f_{0,d}(\frac{1}{b})$ are rational.
For an integer $j$ with $1\le j\le d$, we define integers $p_j,q_j\neq 0$ by $f_{0,j}(1/b)=\frac{p_j}{q_j}$.
By the irrationalty of $f(z)$ and \eqref{aux6}, for any non-negaitive integer $e$, we have
\begin{align}\label{p111}
Q_e(z^{k^e})f(z)-\sum_{j=1}^d a_{e,j,0}(z)P_{e,j}(z^{k^e})=z^{(dp+p+1-J){k^e}}\sum_{j=1}^d a_{e,j,0}(z)G_{e,j}(z^{k^e})\neq 0.
\end{align}
From \eqref{p111}, there exist integers $i_e$ with $1\le i_e\le d$ and $L_e$ such that
\begin{align}\label{p}
ord\; G_{e,i_e}(z^{k^{e}})= L_ek^{e}.
\end{align}
By $(3.5)$, \eqref{p} and \eqref{asss}, for any sufficiently large integer $e$, we get
\begin{align}\label{pq}
0<\frac{1}{b^{L_ek^e}}(1-\frac{e^{\epsilon c_3p^2k^e}}{b^{k^e}-e^{\epsilon k^e}})\le|G_{e,i_e}(\frac{1}{b^{k^e}})|\le\frac{e^{\epsilon c_3p^2k^e}b^{k^e}}{b^{k^e}-e^{\epsilon k^e}}. 
\end{align}
By \eqref{aux11}, we get
\begin{align}\label{pqr}
Q_e(z^{k^e})\left(
    \begin{array}{c}
      f_{0,1}(z) \\
      f_{0,2}(z) \\
      \vdots \\
      f_{0,d}(z)
    \end{array}
  \right)-B_e(z)\left(
    \begin{array}{c}
      P_{e,1}(z^{k^e}) \\
     P_{e,2}(z^{k^e})\\
      \vdots \\
      P_{e,d}(z^{k^e})
    \end{array}
  \right)=B_e(z)z^{(dp+p-J+1)k^e}\left(
    \begin{array}{c}
      G_{e,1}(z^{k^e}) \\
     G_{e,2}(z^{k^e})\\
      \vdots \\
      G_{e,d}(z^{k^e})
    \end{array}
  \right).
\end{align}
We denote the $(i,j)$-component of $B^{-1}_e(\frac{1}{b})$ by $b_{e,i,j}$.
By \eqref{pqr}, we have
\begin{align}\label{r}
(b_{e,i_e,1},\cdots, b_{e,i_e,d})Q_e(\frac{1}{b^{k^e}})\left(
    \begin{array}{c}
      f_{0,1}(\frac{1}{b}) \\
      f_{0,2}(\frac{1}{b}) \\
      \vdots \\
      f_{0,d}(\frac{1}{b})
    \end{array}
  \right)-
      P_{e,i_e}(\frac{1}{b^{k^e}})=\frac{1}{b^{(dp+p-J+1)k^e}}
      G_{e,i_e}(\frac{1}{b^{k^e}}).
  \end{align}
We define a positive integer $C_{1,e}$ by $C_{1,e}:=\min\{ D \;|\;D(b_{e,i_e,1},\cdots, b_{e,i_e,d})\in {\mathbb{Z}}^{d}\}$.
For any non-negaitive integer $e$, we define an integer $I_e$ as follows
\begin{align}\label{s}
I_e:=\prod_{i=1}^{d} q_i C_{1,e} b^{(dp+1-J)k^e}(b_{e,i_e,1},\cdots, b_{e,i_e,d})Q_e(\frac{1}{b^{k^e}} )\left(
    \begin{array}{c}
      f_{0,1}(\frac{1}{b}) \\
      f_{0,2}(\frac{1}{b}) \\
      \vdots \\
      f_{0,d}(\frac{1}{b})
    \end{array}
  \right)-
      \prod_{i=1}^{d} q_i C_{1,e} b^{(dp+1-J)k^e}P_{e,i}(\frac{1}{b^{k^e}}).
  \end{align}
By \eqref{pq}, \eqref{r} and Lemma \ref{invert}, we have
\begin{align}\label{t}
0< |I_e|\le \frac{C_{1,e}\prod_{i=1}^{d} q_i}{b^{pk^e}} \frac{e^{\epsilon c_3p^2k^e}b^{k^e}}{b^{k^e}-e^{\epsilon k^e}}\le
\frac{C_2(\epsilon)d!^e b^{dk^e} e^{\epsilon dk^e}\prod_{i=1}^{d} {q_i}}{{b^{pk^e}}} \frac{e^{\epsilon c_3p^2k^e}b^{k^e}}{b^{k^e}-e^{\epsilon k^e}}.
  \end{align}
By \eqref{t}, Lemma \ref{invert} and $p>d+5$, for any sufficiently large integer $e$, we have
\begin{align}\label{u}
0< |I_e| <1.
  \end{align}
This contradicts that $I_e$ is an integer.
This completes the proof of Theorem \ref{main}.
\hfill
$\square$
\section{The other examples of Theorem \ref{main} by related the digital pattern sequences.}
In this section, we give the examples of Theorem \ref{main} by related the certain digital pattern sequences as follows.
Let $A$ be a set and ${A}^{*}$ be the free monoid generated by $A$.
Let $A$ and $B$ be two finite words on $\mathbb{Z}^*$, and let $AB$ be the concatenation of $A$ and $B$.
 Let $m$ be a positive integer.
 Let $b_i\in \mathbb{Z}$ with $1\le i \le m$, and $A:=b_1b_2\cdots b_m\in {\{b_1,b_2,\cdots b_m\}}^{*}$.
Then, for any element $f\in \mathbb{Z}$,
we define a word $f(A)$ by, $$f(A):=f\cdot b_1 f\cdot  b_2\cdots f\cdot  b_m\in {\{f\cdot b_1,f\cdot b_2,\cdots, f\cdot b_m\}}^{*},$$where, for any $b\in \mathbb{Z}$, $f(b):=f\cdot b\in \mathbb{Z}$ with denoting the multiplication on $\mathbb{Z}$ by the dot symbol. Note that the length of $f(A)$ as well as that of $A$ is $m$. The element $f$ can be regarded as the coding on ${\mathbb{Z}}^*$. (See the definition of coding 9p in [AlS2])
By Theorem $3.2$ in [Mi] (See also example $3.4$ in [Mi]), the following sequences ${(a(n))}_{n\ge0}$ and ${(b(n))}_{n\ge0}$ are generalized $k$-regular.
Let $A_0=a$, $B_0=b$ with $a,b \in \mathbb{Z}$. We define the words $A_{n+1}$ and $B_{n+1}$ of length $2^{n+1}$ recursively as
\begin{align}\label{k-re}
     & A_{n+1}:=A_n B_n,\\ \nonumber
     & B_{n+1}:=A_n f_n(B_n),
\end{align}
where $f_n\in \mathbb{Z}$. For $a=b=1$, we define the ${(a(n))}_{n\ge0}$ by ${(a(n))}_{n\ge0}=\lim_{n \to \infty} A_{n}$.
If $f_n=-1$ $($ $n\ge0$ $)$, then the sequence ${(a(n))}_{n\ge0}$ is known as the Rudin-Shapiro sequence. (See $126$p in [Fo].)
\begin{Rem}{\rm
By th definition of the sequence ${(a(n))}_{n\ge0}$, 
the sequence ${(a(n))}_{n\ge0}$ has the following digital pattern definition. We define the counting functions $d_1(n; 2^y+2^{y+1})$  as
\begin{displaymath}
d_1(n; 2^y+2^{y+1}):=\begin{cases} f_y& \mbox{  $2^y+2^{y+1}$ is}\\ & \mbox{appeared in the base-$2$ representation of $n$  } \\
                      1 & \mbox{Otherwise.}
                                             \end{cases}
\end{displaymath}
The sequence ${(a(n))}_{n\ge0}$ has the following definition
\begin{align}\label{ad}
    a(n) = \prod_{y=0}^\infty d_1(n;2^y+2^{y+1}).
\end{align}}
\end{Rem}
Moreover, for $a=1,b=-1$, we also define the ${(b(n))}_{n\ge0}$ by ${(b(n))}_{n\ge0}=\lim_{n \to \infty} A_{n}$.
Let $f(z):=\sum_{n=0}^\infty a(n)z^n$ and $g(z):=\sum_{n=0}^\infty b(n)z^n$.
From the proof of Theorem $3.2$ in [Mi], $f(z)$ and $g(z)$ has the following representation
\begin{align}\label{rep}
  &\left(\begin{array}{c}
      f(z) \\
      g(z)
      \end{array}
  \right) =\prod_{e=0}^{\underrightarrow{\infty}}\begin{pmatrix}
1+\frac{1+f_e}{2}z^{2^e} & \frac{1-f_e}{2}z^{2^e} \\
1+\frac{-1-f_e}{2}z^{2^e}  &\frac{-1+f_e}{2}z^{2^e}
\end{pmatrix}\left(\begin{array}{c}
      1 \\
      1
      \end{array}
  \right)\nonumber \\
 &=\lim_{e \to \infty}
\begin{pmatrix}
1+\frac{1+f_0}{2}z & \frac{1-f_0}{2}z \\
1+\frac{-1-f_0}{2}z &\frac{-1+f_0}{2}z
\end{pmatrix}
\cdots
\begin{pmatrix}
1+\frac{1+f_e}{2}z^{2^e} & \frac{1-f_e}{2}z^{2^e} \\
1+\frac{-1-f_e}{2}z^{2^e}  &\frac{-1+f_e}{2}z^{2^e}
\end{pmatrix}\left(\begin{array}{c}
      1\\
      1
      \end{array}
  \right).\end{align}
This representation \eqref{rep} show that $f(z)$ and $g(z)$ satisfy the equations \eqref{8} with $A_e(z)=\begin{pmatrix}
1+\frac{1+f_e}{2}z^{} & \frac{1-f_e}{2}z^{} \\
1+\frac{-1-f_e}{2}z^{}  &\frac{-1+f_e}{2}z^{}
\end{pmatrix}$. Now we give the examples of Theorem \ref{main} as follows.
\begin{Prop}\label{exa}
 Notation is same as above. Let $b$ be an integer with $b\ge 2$.
Assume that, for any non-nagative integer $e$, the integer $f_e\neq 1$ is odd.  Moreover, assume that, for any $\epsilon>0$, there exists an integer $N(\epsilon)$ such that, for any $e\ge N(\epsilon)$, 
\begin{align}\label{asy200}
|f_e|\le {e}^{\epsilon k^e}.
\end{align}
Then at least one among the numbers $f(\frac{1}{b})$, $g(\frac{1}{b})$ is transcendental.
\end{Prop}
\begin{proof}
{\rm
By \eqref{k-re}, we have
\begin{align}\label{k-re1}
     & {(a(n))}_{n\ge0}:=A_n B_nA_nf_n(B_n)\cdots.
\end{align}
For any non-negative integer $n$, we define the polynomial $P_n(z)$ as the generating function of $A_nB_n$.
We also define the polynomial $Q_n(z)$ by $Q_n(z)=(1-z^{2^{n+1}})$.
By the definitions of $P_n(z), Q_n(z)$, we have
\begin{align}\label{k-re2}
     \deg P_n(z), \deg Q_n(z)\le 2^{n+1}.
\end{align}
From the definitions of $P_n(z), Q_n(z)$ and \eqref{k-re1},
we have
\begin{align}\label{-re}
   ord (Q_n(z)f(z)-P_n(z) )\ge 2^{n+1}+2^{n}.
\end{align}
By the  $f_n\neq 1$ and the definitions of $P_n(z)$,we have
\begin{align}\label{-re1}
      z^{2^{n+1}}(P_{n+1}(z)-z^{2^{n+1}}P_{n}(z))\neq P_{n+1}(z)-P_{n}(z).
\end{align}
From \eqref{-re1}, we have
\begin{align}\label{-re2}
      P_{n}(z)Q_{n+1}(z)-Q_{n}(z)P_{n+1}(z)\neq 0.
\end{align}
By \eqref{k-re},\eqref{-re},\eqref{-re2} and Lemma \ref{DN}, we have
\begin{align}\label{-re3}
      \mu(f(z))<\infty.
\end{align}
From \eqref{-re3}, $f(z)$ is irrational.
Therefore, by Theorem \ref{main}, at least one among the numbers $f(\frac{1}{b})$, $g(\frac{1}{b})$ is transcendental.}
\end{proof}
By the similar way of Proposition \ref{exa}, one can construct the other examples of Theorem \ref{main} in {\it $k$-recursive sequences}. (See Definition$3.1$ in [Mi]).


\begin{thebibliography}{15}
\bibitem[AmV1]{AmV}M. Amou and K. V\"a\"an\"anen,
\emph{Arithmetical properties of certain infinite products},
J. Number Theory 153 (2015), 283-303.
\bibitem[AmV2]{Am2}M. Amou and K. V\"a\"an\"anen,
\emph{On algebraic independence of a class of infinite products},
J. NumberTheory172(2017)114-–132.
\bibitem[Bec]{Bec}   P. G. Becker,
\emph{k-regular power series and Mahler-type functional equations}, J. Number. Theory. 49 (1994), 269-286.
\bibitem[BelBC]{BeBC} J. P. Bell, Y. Bugeaud and M. Coons,
\emph{Diophantine approximation of Mahler numbers}, Proc. London Math. Soc. (3) 110 (2015) 1157–-1206.
\bibitem[Bu]{Bu}Y. Bugeaud,
\emph{Distribution modulo one and Diophantine approximation}, Cambridge Tracts in Mathematics 193.
\bibitem[CoS]{CoS}M. Coons and L. Spiegelhofer,
\emph{Chapter 2 Number Theoretic Aspects of Regular
Sequences,} Sequences,Groups, and Number Theory Trends in Mathematics (2018) 37--87.
\bibitem[DuN]{DuN}D. Duverney and K. Nishioka,
\emph{An inductive method for proving the transcendence of certain series},
Acta Arith. 110 (4) (2003) 305-330.
 \bibitem[Fo]{Fo} N. P. Fogg,
\emph{Substitutions in Dynamics, Arithmetics and Combinatorics},  
Lecture Notes in Mathematics  Springer; (2002).
\bibitem[Gel]{Gel}Gel'fond, A. O.
\emph{Sur les nombres qui ont des propri\'et\'es additives et multiplicatives donn\'ees},
 Acta Arith. 13 1967/1968 259-265.
\bibitem[Ma]{Ma}D. Masser,
\emph{Auxiliary Polynomials in Number Theory},Cambridge Tracts in Mathematics Book 207.
\bibitem[Mi]{Mi}E. Miyanohara,
\emph{Generalized $k$-regular sequence I: Fundamental properties and examples}, arXiv:1809.09320v2.
\bibitem[Ni]{N}K. Nishioka,
\emph{Mahler Functions and Transcendence},  
LectureNotes in Mathmatics  Springer; (1996).
\bibitem[Ta]{Ta}H. Tachiya,
\emph{Transcendence of certain infinite products}, %
J. Number theorey. 125 (2007), no.1, 182--200.





\end{thebibliography}
\end{document}